\documentclass{article}
\usepackage{amsfonts}
\usepackage{amsmath}
\usepackage{amsmath}
\usepackage{amssymb}

\newtheorem{theorem}{Theorem}

\newtheorem{definition}[theorem]{Definition}

\newtheorem{lemma}[theorem]{Lemma}

\begin{document}

\title{Survival probability for a class of multitype subcritical branching
processes\\
in random environment \thanks{%
This work was supported by the Russian Science Foundation under the grant
17-11-01173}}
\author{Vladimir Vatutin\thanks{%
Novosibirsk State University, Novosibirsk, Russia, e-mail: vatutin@mi-ras.ru}%
, Elena Dyakonova\thanks{%
Novosibirsk State University, Novosibirsk, Russia, e-mail: elena@mi-ras.ru}}
\date{}
\maketitle

\begin{abstract}
We study the asymptotic behaviour of the survival probability of a
multi-type branching processes in random environment. The class of processes
we consider corresponds, in the one-dimensional situation, to the
intermediately subcritical case. We show under rather general assumptions on
the form of the offspring generating functions of particles that the
probability of survival up to generation $n$ of the process initiated at
moment zero by a single particle of any type is of order $\lambda
^{n}n^{-1/2}$ for large $n,$ where $\lambda \in (0,1)$ is a constant
specified by the Lyapunov exponent of the mean matrices of the process.
\end{abstract}



\section{Introduction and main results}

Branching processes in random environment with one type of particles have
been intensively investigated during the last two decades and their
properties are well understood (see, for example, the survey \cite{VZ} and
the recent book by Kersting and Vatutin \cite{KV2017}). The multi-type case
is much less studied and many basic problems such as the asymptotic behavior
of the survival probability, limit theorems for the number of particles in
the process and others are solved under rather heavy conditions, for
example, for the cases when the mean matrices of the reproduction laws of
particles in different generations have a common nonrandom left or right
eigenvector corresponding to their Perron roots, or for some other
relatively narrow classes of mean matrices (see \cite{Dyak18} -- \cite%
{Dyak07}, \cite{Vat2013b}).

This paper supplements some recent results (see \cite{Dy2013}, \cite{LPP2016}%
, \cite{VD2017}, \cite{VatWacht2018}) describing the asymptotic behavior of
the survival probabilities of the critical and subcritical multitype
branching processes evolving in random environment.

To formulate our main result we need some notation for $p$-dimensional
vectors and $p\times p$ matrices. We usually make no difference in notation
for row and column vectors. As we hope it will be clear from the context
which form is selected in each case. Besides we write

$\mathbf{e}_{j}$, $j=1,\dots ,p$, for a vector whose $j$-th component is
equal to $1$ and the others are zeros;

$\mathbf{0}=(0,\dots ,0)$, $\mathbf{1}=(1,\dots ,1)$ for zero and unit $p$%
-dimensional vectors.

The norm and scalar product of vectors $\mathbf{x}=(x_{1},\dots ,x_{p})$ and
$\mathbf{y}=(y_{1},\dots ,y_{p})$ are denoted as
\begin{equation*}
|\mathbf{x|}=\sum_{i=1}^{p}|x_{i}|,\qquad (\mathbf{x},\mathbf{y}%
)=\sum_{i=1}^{p}x_{i}y_{i}.\qquad
\end{equation*}%
We also use the notation $\mathbf{x}^{\mathbf{y}}=\prod%
\limits_{i=1}^{p}x_{i}^{y_{i}}$ and~define the norm of a matrix $\mathbf{m}%
=(m(i,j))_{i,j=1}^{p}$ as
\begin{equation*}
|\mathbf{m}|=\sum_{i=1}^{p}\sum_{j=1}^{p}|m(i,j)|.
\end{equation*}

Let $\mathcal{P}(\mathbb{N}_{0}^{p})$~be the space of all probability
measures on the set $\mathbb{N}_{0}^{p}$ of $p$-dimensional vectors with
nonnegative integer-valued components. For a measure $f\in \mathcal{P}(%
\mathbb{N}_{0}^{p})$ we denote by $f[\mathbf{z}]$ the mass assigning by the
measure to the point $\mathbf{z=}(z_{1},\dots ,z_{p})\in \mathbb{N}_{0}^{p}$%
. The function
\begin{equation*}
f\mathbf{(s):=}\sum_{\mathbf{z}\in \mathbb{N}_{0}^{p}}f[\mathbf{z}]\mathbf{s}%
^{\mathbf{z}},\qquad \mathbf{s=}(s_{1},\dots ,s_{p})\in \lbrack 0,1]^{p},
\end{equation*}%
is the generating function for the distribution (measure) $f$. It will be
convenient to denote (by taking some liberty) the distribution (measure) and
the corresponding generating function by one and the same symbol $f$. We
also need $p$-dimensional vectors
\begin{equation*}
\mathbf{f}=(f^{(1)},\dots ,f^{(p)})\in \mathcal{P}(\mathbb{N}_{0}^{p})\times
\dots \times \mathcal{P}(\mathbb{N}_{0}^{p})=:\mathcal{P}^{p}(\mathbb{N}%
_{0}^{p}),
\end{equation*}%
whose components are probability measures $f^{(i)}\in \mathcal{P}(\mathbb{N}%
_{0}^{p})$, $i=1,\dots ,p$. In what follows it will be sometimes convenient
to call vectors $\mathbf{f}\in \mathcal{P}^{p}(\mathbb{N}_{0}^{p})$ simply
as probability measures and the corresponding vectors $\mathbf{f(s)}$ of
generating functions as generating functions.

\begin{definition}
\label{d1} A sequence $v=(\mathbf{f}_{1},\mathbf{f}_{2},\dots )$ of
probability measures on $(\mathbb{N}_{0}^{p})^{p}$ is called a \textit{%
varying environment}.
\end{definition}

\begin{definition}
\label{d2} Let $v=(\mathbf{f}_{n},n\geq 1)$~be a varying environment. A
stochastic process $\bigl\{\mathbf{Z}_{n}=(Z_{n}(1),\dots ,Z_{n}(p)),\,n\geq
0\bigr\}$ with values in the space $\mathbb{N}_{0}^{p}$ is called a \textit{%
branching process in the environment} $v$, if, for any $\mathbf{z}\in
\mathbb{N}_{0}^{p}$ and $n\geq 1$
\begin{equation*}
\mathbf{P}(\mathbf{Z}_{n}=\mathbf{z}\mid \mathbf{Z}_{0},\dots ,\mathbf{Z}%
_{n-1})=(\mathbf{f}_{n}^{\mathbf{Z}_{n-1}})[\mathbf{z}].
\end{equation*}
\end{definition}

In the sequel the symbol $\mathbf{P}_{\mathbf{z},v}(\,\cdot \,)$ will
correspond to the distribution of the process in the varying environment $v$
under the initial value $\mathbf{Z}_{0}=\nobreak\mathbf{z}$.

We now introduce the notion of a multitype branching process in random
environment specified on the corresponding probability space $(\Omega ,%
\mathcal{F},\mathbb{P})$. Define on the set $\mathcal{P}^{p}(\mathbb{N}%
_{0}^{p})$ of probability measures the metric of total variation $d_{\mathrm{%
TV}}$ by the formula
\begin{equation*}
d_{\mathrm{TV}}(\mathbf{f,g})=\frac{1}{2p}\sum_{\mathbf{z}\in \mathbb{N}%
_{0}^{p}}|\mathbf{f}[\mathbf{z}]-\mathbf{g}[\mathbf{z}]|,\qquad \mathbf{f,g}%
\in \mathcal{P}^{p}(\mathbb{N}_{0}^{p}),
\end{equation*}%
and supply $\mathcal{P}^{p}(\mathbb{N}_{0}^{p})$ with the Borel $\sigma $%
-algebra generated by $d_{\mathrm{TV}}$.

We consider random probability measures $\mathbf{F=}(F^{(1)},\dots ,F^{(p)})$
being random vectors with values in the space $\mathcal{P}^{p}(\mathbb{N}%
_{0}^{p})$, whose components are specified by the probability generating
functions in $p$ variables:
\begin{equation*}
F^{(i)}(\mathbf{s}):=\sum_{\mathbf{z}\in \mathbb{N}_{0}^{p}}F^{(i)}[\mathbf{z%
}]\mathbf{s}^{\mathbf{z}},\qquad i=1,\dots ,p.
\end{equation*}

\begin{definition}
\label{d3} A sequence $\mathcal{V}$ $=\{\mathbf{F}_{1},\mathbf{F}_{2},\dots
\}$ of random measures is called a \textit{random environment}.
\end{definition}

We say that the random environment $\mathcal{V}$ is generated by a sequence
of independent identically distributed random variables if the random
measures $\mathbf{F}_{1},\mathbf{F}_{2},\dots $ are independent copies of a
random probability measure $\mathbf{F}$ with values in~$\mathcal{P}^{p}(%
\mathbb{N}_{0}^{p})$. In this paper we deal with such an environment only.

In what follows the symbols $\mathbb{P}$ and $\mathbb{E}$ denote probability
and expectation for a branching process in a random environment in contrast
to the symbols $\mathbf{P}$ and $\mathbf{E}$ applied in the case of a
branching process in a varying environment.

\begin{definition}
\label{d4} Let $\mathcal{V}$~be a random environment. A stochastic process
\begin{equation*}
\mathcal{Z}=\left\{ \mathbf{Z}_{n}=(Z_{n}(1),\dots ,Z_{n}(p)),n\geq 0\right\}
\end{equation*}%
with values in~$\mathbb{N}_{0}^{p}$ is called a $p$-\textit{type branching
process} in the random environment $\mathcal{V}$, if, for all $\mathbf{z},%
\mathbf{z}_{1},\dots ,\mathbf{z}_{k}\in \mathbb{N}_{0}^{p}$ and any fixed
environment $v$
\begin{align*}
& \mathbb{P}(\mathbf{Z}_{1}=\mathbf{z}_{1},\dots ,\mathbf{Z}_{k}=\mathbf{z}%
_{k}\mid \mathbf{Z}_{0}=\mathbf{z};\mathcal{V}=v) \\
& \qquad =\mathbf{P}_{\mathbf{z},v}(\mathbf{Z}_{1}=\mathbf{z}_{1},\dots ,%
\mathbf{Z}_{k}=\mathbf{z}_{k})\quad \mathbb{P}\text{-a.s.}
\end{align*}
\end{definition}

We use below the uppercase letters to denote variables or functions if we
deal with a random environment, and the lowercase letters to denote the
corresponding variables or functions if we deal with a fixed environment.
For instance, the (random) distribution law of particles of the $(n-1)$th
generation will be specified by a tuple $\mathbf{F}_{n}=(F_{n}^{(1)},\dots
,F_{n}^{(p)})$ of (random) probability generating functions in $p$
variables. Similarly, we denote by
\begin{equation*}
\mathbf{M}_{n}:=(M_{n}(i,j))_{i,j=1}^{p}=\biggl(\frac{\partial F_{n}^{(i)}}{%
\partial s_{j}}(\mathbf{1})\biggr)_{i,j=1}^{p}
\end{equation*}%
the mean matrix corresponding to the probability generating function $%
\mathbf{F}_{n}$, and so on. Clearly, the random matrices $\mathbf{M}_{n}$, $%
n\geq 1$, as well as the matrix
\begin{equation*}
\mathbf{M}=(M(i,j))_{i,j=1}^{p}:=\biggl(\frac{\partial F^{(i)}}{\partial
s_{j}}(\mathbf{1})\biggr)_{i,j=1}^{p}
\end{equation*}%
are independent and identically distributed under our conditions.

We define the cone
\begin{equation*}
\mathcal{C}:\mathcal{=}\{\mathbf{x}=(x_{1},\dots ,x_{p})\in \mathbb{R}%
^{p}\colon x_{i}\geq 0\text{ for all }i=1,\dots ,p\},
\end{equation*}%
the sphere
\begin{equation*}
\mathbb{S}^{p-1}:=\{\mathbf{x}\colon \mathbf{x}\in \mathbb{R}^{p},|\mathbf{x}%
|=1\},
\end{equation*}%
and the space $\mathbb{X}:=\mathcal{C}\cap \mathbb{S}^{p-1}$. In the sequel
we need to consider the linear semi-group $S^{+}$ of $p\times p$ matrices
with nonnegative elements each whose row and column includes at least one
positive element. For a vector $\mathbf{x}\in \mathbb{X}$ and a matrix $%
\mathbf{m}\in S^{+}$ we specify the projective actions
\begin{equation*}
\mathbf{x}\cdot \mathbf{m}:=\frac{\mathbf{xm}}{|\mathbf{xm}|},\quad \mathbf{m%
}\cdot \mathbf{x}:=\frac{\mathbf{mx}}{|\mathbf{mx}|}
\end{equation*}%
and define a function $\rho $ on the product space $\mathbb{X}\times
S^{+}=\{(\mathbf{x},\mathbf{m})\}$ by setting
\begin{equation*}
\rho (\mathbf{x},\mathbf{m}):=\log |\mathbf{xm}|.
\end{equation*}%
The function meets the so-called cocycle property meaning that for a vector $%
\mathbf{x}\in \mathbb{X}$ and matrices $\mathbf{m}_{1},\mathbf{m}_{2}\in
S^{+}$
\begin{equation*}
\rho (\mathbf{x},\mathbf{m}_{1}\mathbf{m}_{2})=\rho (\mathbf{x\cdot m}_{1},%
\mathbf{m}_{2})+\rho (\mathbf{x},\mathbf{m}_{1}).
\end{equation*}

The measure $\mathbb{P}$, generated by a branching process in random
environment (BPRE) with $p$ types of particles, specifies the corresponding
probability measure on the Borel $\sigma $-algebra of the semi-group $S^{+}.$
We agree to denote this measure as $\mathbb{P}$ as well, i.e., for a Borel
subset $\mathcal{A}\subseteq \nobreak S^{+}$ we set
\begin{equation*}
\mathbb{P}(\mathbf{M}\in \mathcal{A}):=\mathbb{P}(\mathbf{f}\colon \mathbf{M}%
=\mathbf{M}(\mathbf{f})\in \mathcal{A}).
\end{equation*}

Keeping in mind this agreement we introduce a number of assumptions to be
valid throughout the paper. These assumptions are simplified versions of the
conditions introduced in \cite{PH17} and concern only properties of the
restriction of $\mathbb{P}$ to the semi-group~$S^{+}$. \

\begin{itemize}
\item \textbf{Condition }$\mathbf{H1}$. The set $\Theta :=\left\{ \theta >0:%
\mathbb{E}\left[ \left\vert \mathbf{M}\right\vert ^{\theta }\right] <\infty
\right\} $ is nonempty.

\item \textbf{Condition }$\mathbf{H2}$. There exists a positive number $%
\Delta >1$ such that
\begin{equation*}
1\leq \frac{\max_{i,j}M(i,j)}{\min_{i,j}M(i,j)}\leq \Delta .
\end{equation*}

\item \textbf{Condition }$\mathbf{H3}$. There exists $\delta >0$ such that%
\begin{equation*}
\inf_{\mathbf{x}\in \mathbb{X}}\mathbb{P}\left( \mathbf{M}:\log \left\vert
\mathbf{Mx}\right\vert >\delta \right) >0.
\end{equation*}
\end{itemize}

Along with random matrices $\mathbf{M}_{n}$ and $\mathbf{M}$ we introduce
the random Hessian matrices
\begin{equation}
\mathbf{B}^{(i)}:=\biggl(\frac{\partial ^{2}F^{(i)}}{\partial
s_{k}\,\partial s_{l}}(\mathbf{1})\biggr)_{k,l=1}^{p},\quad \mathbf{B}%
_{n}^{(i)}:=\biggl(\frac{\partial ^{2}F_{n}^{(i)}}{\partial s_{k}\,\partial
s_{l}}(\mathbf{1})\biggr)_{k,l=1}^{p},  \notag
\end{equation}%
and set
\begin{equation*}
\mathcal{B}:=\sum_{i=1}^{p}|\mathbf{B}^{(i)}|,\quad \mathcal{T}:=\frac{%
\mathcal{B}}{|\mathbf{M}|^{2}},\quad \mathcal{B}_{n}:=\sum_{i=1}^{p}|\mathbf{%
B}_{n}^{(i)}|,\quad \mathcal{T}_{n}:=\frac{\mathcal{B}_{n}}{|\mathbf{M}%
_{n}|^{2}}.
\end{equation*}

Thus, $\mathcal{T}_{n}$ are independent probabilistic copies of $\mathcal{T}$%
. We shall impose, along with Conditions $\mathbf{H1}-\mathbf{H3}$ the
following restriction on the distribution of $\mathcal{T}$.

\begin{itemize}
\item \textbf{Condition }$\mathbf{H4}$. There exists an $\varepsilon >0$
such that
\begin{equation*}
\mathbb{E}\left[ \left\vert \mathbf{M}\right\vert |\log \mathcal{T}%
|^{1+\varepsilon }\right] <\infty .
\end{equation*}
\end{itemize}

Using the standard subadditivity arguments, one can easily infer that for
every $\theta \in \Theta $ the limit
\begin{equation*}
\lambda \left( \theta \right) :=\lim_{n\rightarrow \infty }\left( \mathbb{E}%
\left[ \left\vert \mathbf{M}_{n}\cdot \cdot \cdot \mathbf{M}_{1}\right\vert
^{\theta }\right] \right) ^{1/n}<\infty
\end{equation*}%
is well defined. This function is an analog of \ the moment generating
function for the associated random walk in the case of single-type BPRE's.

Set
\begin{equation*}
\Lambda (\theta ):=\log \lambda (\theta ),\quad \theta \in \Theta .
\end{equation*}
Here is our main result.

\begin{theorem}
\label{thm:nonextinct} Assume that Conditions $\mathbf{H1}-\mathbf{H4}$ are
valid, the point $\theta =1$ belongs to the interior of the set $\Theta $
and, in addition, $\Lambda ^{\prime }(0)<0$ and $\Lambda ^{\prime }(1)=0$.
Then there exist positive constants $C^{-}$ and $C^{+}$ such that, for all $%
i=1,...,p$ and all $n\geq 1$
\begin{equation}
\frac{C^{-}}{\sqrt{n}}\,\lambda ^{n}(1)\leq \mathbb{P}\left( |\mathbf{Z}%
_{n}|>0\Big|\mathbf{Z}_{0}=\mathbf{e}_{i}\right) \leq \frac{C^{+}}{\sqrt{n}}%
\,\lambda ^{n}(1).  \label{P_asymp}
\end{equation}
\end{theorem}

Dyakonova \cite{Dy2013} has proved a statement more precise than (\ref%
{P_asymp}) under stronger restrictions. Namely, she has shown that if all
possible realisations of $\mathbf{M}$ have a common deterministic left
eigen-vector $\mathbf{v}$ corresponding to the Perron root $\chi (\mathbf{M})
$ of $\mathbf{M}$ and some other technical conditions are valid then there
exists a vector $\mathbf{C=}\left( C_{1},...,C_{p}\right) $ with strictly
positive components such that,%
\begin{equation*}
\mathbb{P}\left( |\mathbf{Z}_{n}|>0\Big|\mathbf{Z}_{0}=\mathbf{e}_{i}\right)
\sim \frac{C_{i}}{\sqrt{n}}\,\lambda ^{n}(1),\quad n\rightarrow \infty .
\end{equation*}

Note that the assumption $\Lambda ^{\prime }(1)=0$ reduces in this special
case to the condition $\mathbb{E}[\chi (\mathbf{M})\log \chi (\mathbf{M})]=0.
$ In the single-type case the last condition corresponds  to the so-called
intermediately subcritical BPRE's (see, for instance, \cite{ABKV2014} or
\cite{KV2017}, chapter 8).

\section{Auxiliary results}

Denote by $\mathcal{C}\left( \mathbb{X}\right) $ the set of all continuous
functions on $\mathbb{X}$. For $\theta \in \Theta ,$ $g\in \mathcal{C}\left(
\mathbb{X}\right) ,$ and $\mathbf{x\in }\mathbb{X}$ define the transition
operators
\begin{equation*}
P_{\theta }g(\mathbf{x}):=\mathbb{E}\left[ \left\vert \mathbf{Mx}\right\vert
^{\theta }g\left( \mathbf{M\cdot x}\right) \right]
\end{equation*}%
and
\begin{equation*}
P_{\theta }^{\ast }g(\mathbf{x}):=\mathbb{E}\left[ \left\vert \mathbf{M}^{T}%
\mathbf{x}\right\vert ^{\theta }g\left( \mathbf{M}^{T}\mathbf{\cdot x}%
\right) \right] ,
\end{equation*}%
where $\mathbf{M}^{T}$ is the matrix transposed to $\mathbf{M}$.

If Conditions $\mathbf{H1}-\mathbf{H3}$ hold, then, according to Proposition
3.1 in~\cite{BDGM2014}, $\lambda (\theta )$ is the spectral radius of $%
P_{\theta }$ and $P_{\theta }^{\ast }$ and there exist unique strictly
positive functions $r_{\theta },r_{\theta }^{\ast }\in \mathcal{C}\left(
\mathbb{X}\right) $ and unique probability measures $l_{\theta }$ and $%
l_{\theta }^{\ast }$ subject to the scalings
\begin{equation*}
\int_{\mathbb{X}}r_{\theta }(\mathbf{x})dl_{\theta }(\mathbf{x})=1,\quad
\int_{\mathbb{X}}r_{\theta }^{\ast }(\mathbf{x})dl_{\theta }^{\ast }(\mathbf{%
x})=1
\end{equation*}%
and possessing the properties
\begin{equation}
l_{\theta }P_{\theta }=\lambda \left( \theta \right) l_{\theta },\quad
P_{\theta }r_{\theta }=\lambda \left( \theta \right) r_{\theta }.
\label{dop}
\end{equation}%
\begin{equation*}
l_{\theta }^{\ast }P_{\theta }^{\ast }=\lambda \left( \theta \right)
l_{\theta }^{\ast },\quad P_{\theta }^{\ast }r_{\theta }^{\ast }=\lambda
\left( \theta \right) r_{\theta }^{\ast }.
\end{equation*}

Following \cite{CM2016}, we introduce the functions
\begin{equation}
p_{n}^{\theta }\left( \mathbf{x},\mathbf{m}\right) :=\frac{\left\vert
\mathbf{mx}\right\vert ^{\theta }}{\lambda ^{n}\left( \theta \right) }\frac{%
r_{\theta }\left( \mathbf{m\cdot x}\right) }{r_{\theta }\left( \mathbf{x}%
\right) },\quad \mathbf{x}\in \mathbb{X}\text{.}  \label{p-def}
\end{equation}

It is easy to see that, for $n\geq 1$, $\mathbf{x}\in \mathbb{X}$ and $%
\mathbf{m\in }S^{+}$
\begin{equation}
\mathbb{E}\left[ p_{n+1}^{\theta }\left( \mathbf{x},\mathbf{Mm}\right) %
\right] =p_{n}^{\theta }\left( \mathbf{x},\mathbf{m}\right)
\label{consistency1}
\end{equation}%
and, in view of (\ref{dop})
\begin{equation}
\mathbb{E}\left[ p_{n}^{\theta }\left( \mathbf{x},\mathbf{L}_{n,1}\right) %
\right] =1.  \label{total_mass}
\end{equation}%
For each $n\geq 1$ let $\mathcal{F}_{n}$ be the $\sigma $-algebra generated
by random elements $\mathbf{Z}_{1},\mathbf{Z}_{2},\ldots ,\mathbf{Z}_{n}$
and $\mathbf{F}_{1},\mathbf{F}_{2},\ldots ,\mathbf{F}_{n}$. It follows from %
\eqref{total_mass} that
\begin{equation*}
\mathbb{P}_{n}^{\theta }(A):=\mathbb{E}\left[ p_{n}^{\theta }\left( \mathbf{x%
},\mathbf{L}_{n,1}\right) \mathbb{I}_{A}\right]
\end{equation*}%
is a probability measure on $\mathcal{F}_{n}$ (here $\mathbb{I}_{A}$ is the
indicator of the event $A$). Furthermore, \eqref{consistency1} implies that
the sequence of measures $\left\{ \mathbb{P}_{n}^{\theta },n\geq 1\right\} $
is consistent and can be extended to a probability measure $\mathbb{P}%
^{\theta }$ on our original probability space $(\Omega ,\mathcal{F})$.
Denote by $\mathbb{E}^{\theta }\left[ \cdot \right] $ the expectation taken
with respect to this measure.

Now we take $\theta =1$ and introduce a homogeneous Markov chain $\left\{
\mathbf{X}_{n},n\geq 0\right\} $ with values in $\mathbb{X}$, where
\begin{equation*}
\mathbf{X}_{0}:=\mathbf{x}\in \mathbb{X}\text{ and }\mathbf{X}_{n}:=\mathbf{%
x\cdot }\left( \mathbf{L}_{n,1}\right) ^{T},\;n\geq 1.
\end{equation*}%
Observe that $\left\vert \mathbf{x\cdot }\left( \mathbf{L}_{n,1}\right)
^{T}\right\vert >0$ by Condition $\mathbf{H2}.$ Since the matrices $\mathbf{M%
}_{n\text{ }}$ are i.i.d. with respect to the measure $\mathbb{P}^{1}\mathbf{%
,}$ the transition probabilities of the chain are specified, for any vector $%
\mathbf{x}\in \mathbb{X}$ and any Borel function $\phi :\mathbb{X}%
\rightarrow \mathbb{R}$ by the relation
\begin{equation*}
Q\phi (\mathbf{x}):=\int_{S^{+}}\phi \left( \mathbf{x}\cdot \mathbf{m}%
^{T}\right) \mathbb{P}^{1}\left( d\mathbf{m}\right) .
\end{equation*}

We fix a vector $\mathbf{x}\in \mathbb{X}$, a number $a<0$ and introduce a
sequence $\left\{ S_{n},n\geq 0\right\} $ by the equalities
\begin{equation*}
\quad S_{0}=a,\quad S_{n}=S_{0}+\log \left\vert \mathbf{xM}_{1}^{T}\mathbf{M}%
_{2}^{T}\cdot \cdot \cdot \mathbf{M}_{n}^{T}\right\vert ,\quad n\geq 1.
\end{equation*}%
Denote $\mathbb{P}_{\mathbf{x},a}^{1}$ the conditional measure, generated by
the measure $\mathbb{P}^{1},$ and $\mathbb{E}_{\mathbf{x},a}^{1}\left[ \cdot %
\right] $ the corresponding conditional expectation given the event $\left\{
\mathbf{X}_{0}=\mathbf{x},S_{0}=a\right\} $.

Let%
\begin{equation*}
\mu :=\min \left\{ n\geq 1:S_{n}\geq 0\right\}
\end{equation*}%
be the first moment when the sequence $\left\{ S_{n},n\geq 1\right\} $
enters the set $[0,\infty ).$

Modifying in a natural way the arguments used in \cite{PH17} or in Appendix
to~\cite{LPP2016} one can conclude that given the conditions of
Theorem \ref{thm:nonextinct} the function $h:$ $\mathbb{X\times (-\infty }%
,0)\rightarrow \mathbb{[}0,\infty )$, specified by the relation
\begin{equation*}
h(\mathbf{x},a):=\lim_{n\rightarrow \infty }\mathbb{E}_{\mathbf{x},a}\left[
-S_{n};\mu >n\right] ,
\end{equation*}%
possesses the property%
\begin{equation}
\mathbb{E}_{\mathbf{x},a}\left[ h(\mathbf{X}_{1},S_{1});\mu >1\right] =h(%
\mathbf{x},a).  \label{HarmonH}
\end{equation}

We need the following upper and lower estimates for $h(\mathbf{x},a)$ which
are reformulations  to our setting the respective results from \cite{PH17}.

\begin{lemma}
\label{L_abbelow} (compare with Theorem 1.1. in \cite{PH17}) Under
Conditions $\mathbf{H1}-\mathbf{H3}$, there exist constants $R>0$ and $%
0<C<\infty $ such that, for all $\left( \mathbf{x},a\right) \in \mathbb{X}%
\times \mathbb{(-\infty },0\mathbb{)}$
\begin{equation}
\max \{C^{-1},\left\vert a\right\vert -R\}<h(\mathbf{x},a)\leq
C(1+\left\vert a\right\vert )  \label{HarmonAbove}
\end{equation}%
and
\begin{equation}
1+\left\vert a\right\vert \leq (R+1)(1+h(\mathbf{x},a)).
\label{HarmonAbove2}
\end{equation}
\end{lemma}

The next result is a restatement of a part of Theorem 1.2 from \cite{PH17}.

\begin{lemma}
\label{L_impor} Let Conditions $\mathbf{H1}-\mathbf{H3}$ be valid. Then, for
any pair $\left( \mathbf{x},a\right) \in \mathbb{X}\times \mathbb{(-\infty }%
,0\mathbb{)}$ as $n\rightarrow \infty $
\begin{equation*}
\mathbb{P}_{\mathbf{x},a}^{1}\left( \mu >n\right) \sim \frac{2}{\sigma \sqrt{%
2\pi n}}\,h(\mathbf{x},a),
\end{equation*}%
where $\sigma \in (0,\infty )$ is a constant. Moreover, there exists a
constant $C>0$ such that, for any pair $\left( \mathbf{x},a\right) \in
\mathbb{X}\times \mathbb{(-\infty },0\mathbb{)}$
\begin{equation}
\mathbb{P}_{\mathbf{x},a}\left( \mu >n\right) \leq \frac{C\left(
1+\left\vert a\right\vert \right) }{\sqrt{n}}  \label{Mu_above}
\end{equation}%
for all $n\geq 1$.
\end{lemma}

Recall that $\mathcal{F}_{n},n\geq 1,$ is the $\sigma -$algebra generated by
the random variables $\mathbf{Z}_{0},\mathbf{Z}_{1},...,\mathbf{Z}_{n},%
\mathbf{F}_{1},\mathbf{F}_{2},...,\mathbf{F}_{n}.$

We introduce a new measure $\widehat{\mathbb{P}}^{1}$ on the flow of $\sigma
-$algebras $\{\mathcal{F}_{n},\,n\geq 1\}$ by setting%
\begin{equation*}
\widehat{\mathbb{E}}_{\mathbf{x},a}^{1}\left[ Y_{n}\right] :=\frac{1}{h(%
\mathbf{x},a)}\mathbb{E}_{\mathbf{x},a}^{1}\left[ Y_{n}h(\mathbf{X}%
_{n},S_{n});\mu >n\right]
\end{equation*}
for any $\left( \mathbf{x},a\right) \in \mathbb{X}\times \mathbb{(-\infty },0%
\mathbb{)}$ and nonnegative random variable $Y_{n}$ measurable with respect
to the $\sigma -$algebra $\mathcal{F}_{n}$.

It follows from (\ref{HarmonH}) and the Markov property that the respective
measure$\widehat{\text{ }\mathbb{P}}_{\mathbf{x},a}^{1}$  is well defined
(compare with the similar definition in \cite{PH17}).

\begin{lemma}
\label{L_basic} Let Conditions $\mathbf{H1}-\mathbf{H3}$ be valid and $Y_{k}$
be a random variable measurable with respect to the $\sigma -$algebra $%
\mathcal{F}_{k},$ $k\geq 1$. Then, for any pair $\left( \mathbf{x},a\right)
\in \mathbb{X}\times (-\infty ,0)\ $%
\begin{equation}
\lim_{n\rightarrow \infty }\mathbb{E}_{\mathbf{x},a}^{1}\left[ Y_{k}|\mu >n%
\right] =\widehat{\mathbb{E}}_{\mathbf{x},a}^{1}\left[ Y_{k}\right] .
\label{FixedK}
\end{equation}%
Moreover, if $Y_{1},Y_{2},...$ is a sequence of uniformly bounded random
variables adopted to the filtration $\left\{ \mathcal{F}_{n},n\geq 1\right\}
$ and  converging $\widehat{\mathbb{P}}_{\mathbf{x},a}^{1}$ - a.s. as $%
n\rightarrow \infty $ to a random variable $Y_{\infty },$ then
\begin{equation}
\lim_{n\rightarrow \infty }\mathbb{E}_{\mathbf{x},a}^{1}\left[ Y_{n}|\mu >n%
\right] =\widehat{\mathbb{E}}_{\mathbf{x},a}^{1}\left[ Y_{\infty }\right] .
\label{Growing_n}
\end{equation}
\end{lemma}

\textbf{Proof}. We follow with minor changes the line of proving lemma 2.5
in \cite{4h} (see also Lemma 5.2 in \cite{KV2017}). Let%
\begin{equation*}
\mathfrak{m}_{x,a}(n):=\mathbb{P}_{\mathbf{x},a}^{1}\left( \mu >n\right) .
\end{equation*}%
Clearly,%
\begin{eqnarray*}
\mathbb{E}_{\mathbf{x},a}^{1}\left[ Y_{k}|\mu >n\right] &=&\frac{1}{\mathbb{P%
}_{\mathbf{x},a}^{1}\left( \mu >n\right) }\mathbb{E}_{\mathbf{x},a}^{1}\left[
Y_{k};\mu >n\right] \\
&=&\frac{1}{\mathbb{P}_{\mathbf{x},a}^{1}\left( \mu >n\right) }\mathbb{E}_{%
\mathbf{x},a}^{1}\left[ Y_{k}\mathfrak{m}_{\mathbf{X}_{k},S_{k}}\left(
n-k\right) ;\mu >k\right] .
\end{eqnarray*}%
In view of Lemma \ref{L_impor}
\begin{equation*}
\lim_{n\rightarrow \infty }\frac{\mathfrak{m}_{\mathbf{X}_{k},S_{k}}\left(
n-k\right) }{\mathbb{P}_{\mathbf{x},a}^{1}\left( \mu >n\right) }=\frac{h(%
\mathbf{X}_{k},S_{k})}{h(\mathbf{x},a)}\text{ }\quad \mathbb{P}_{\mathbf{x}%
,a}^{1}\text{-a.s.}
\end{equation*}%
and there exists a constant $C>0$ such that%
\begin{equation*}
\frac{\mathfrak{m}_{\mathbf{X}_{k},S_{k}}\left( n-k\right) }{\mathbb{P}_{%
\mathbf{x},a}^{1}\left( \mu >n\right) }\leq C\frac{h(\mathbf{X}_{k},S_{k})}{%
h(\mathbf{x},a)}.
\end{equation*}%
The estimate%
\begin{equation*}
\mathbb{E}_{\mathbf{x},a}^{1}\left[ Y_{k}\frac{h(\mathbf{X}_{k},S_{k})}{h(%
\mathbf{x},a)};\mu >k\right] =\widehat{\mathbb{E}}_{\mathbf{x},a}^{1}\left[
Y_{k}\right] <\infty
\end{equation*}%
allows us to apply the dominated convergence theorem to get%
\begin{eqnarray*}
\lim_{n\rightarrow \infty }\mathbb{E}_{\mathbf{x},a}^{1}\left[ Y_{k}|\mu >n%
\right] &=&\mathbb{E}_{\mathbf{x},a}^{1}\left[ Y_{k}\lim_{n\rightarrow
\infty }\frac{\mathfrak{m}_{\mathbf{X}_{k},S_{k}}\left( n-k\right) }{\mathbb{%
P}_{\mathbf{x},a}^{1}\left( \mu >n\right) };\mu >k\right] \\
&=&\frac{1}{h(\mathbf{x},a)}\mathbb{E}_{\mathbf{x},a}^{1}\left[ Y_{k}h(%
\mathbf{X}_{k},S_{k});\mu >k\right] =\widehat{\mathbb{E}}_{\mathbf{x},a}^{1}%
\left[ Y_{k}\right] ,
\end{eqnarray*}%
proving (\ref{FixedK}).

To check the validity of (\ref{Growing_n}) fix $\gamma >1,$ assume for
simplicity that $\left\vert Y_{n}\right\vert \leq 1$ for all $n\geq 1,$ and
observe that in view of Lemma~\ref{L_impor}
\begin{eqnarray*}
\lim_{\gamma \downarrow 1}\lim_{n\rightarrow \infty }\sup \frac{\left\vert
\mathbb{E}_{\mathbf{x},a}^{1}\left[ Y_{n};\mu >n,\mu \leq \gamma n\right]
\right\vert }{\mathbb{P}_{\mathbf{x},a}^{1}\left( \mu >n\right) } &\leq
&\lim_{\gamma \downarrow 1}\lim_{n\rightarrow \infty }\sup \frac{\mathbb{P}_{%
\mathbf{x},a}^{1}\left( \mu >n,\mu \leq \gamma n\right) }{\mathbb{P}_{%
\mathbf{x},a}^{1}\left( \mu >n\right) } \\
&\leq &\lim_{\gamma \downarrow 1}\left( 1-\gamma ^{-1/2}\right) =0.
\end{eqnarray*}%
Further we write for sufficiently large $n$
\begin{eqnarray*}
&&\Delta _{k,n}\left( \gamma \right) :=\left\vert \mathbb{E}_{\mathbf{x}%
,a}^{1}\left[ Y_{n}-Y_{k}|\mu >\gamma n\right] \right\vert =\mathbb{E}_{%
\mathbf{x},a}^{1}\left[ \left\vert Y_{n}-Y_{k}\right\vert \frac{\mathfrak{m}%
_{\mathbf{X}_{n},S_{n}}\left( \left( \gamma -1\right) n\right) }{\mathbb{P}_{%
\mathbf{x},a}^{1}\left( \mu >\gamma n\right) };\mu >n\right]  \\
&&\qquad \qquad \qquad \qquad \qquad \quad \leq 2C\sqrt{\frac{\gamma }{%
\gamma -1}}\frac{1}{h(\mathbf{x},a)}\mathbb{E}_{\mathbf{x},a}^{1}\left[
\left\vert Y_{n}-Y_{k}\right\vert h(\mathbf{X}_{n},S_{n});\mu >n\right]  \\
&&\qquad \qquad \qquad \qquad \qquad \qquad \qquad \qquad \qquad =2C\sqrt{%
\frac{\gamma }{\gamma -1}}\widehat{\mathbb{E}}_{\mathbf{x},a}^{1}\left[
\left\vert Y_{n}-Y_{k}\right\vert \right] .
\end{eqnarray*}%
Since $Y_{m}\rightarrow Y_{\infty }$ $\widehat{\mathbb{P}}_{\mathbf{x},a}^{1}
$ - a.s. as $m\rightarrow \infty $ by the conditions of the lemma, letting
first $n$ to inifinity and than $k$ to infinity $\Delta _{k,n}\left( \gamma
\right) $ vanishes  for any fixed $\gamma >1$. This fact and the first part
of the lemma show that
\begin{equation*}
\lim_{n\rightarrow \infty }\mathbb{E}_{\mathbf{x},a}^{1}\left[ Y_{n}|\mu
>\gamma n\right] =\lim_{k\rightarrow \infty }\lim_{n\rightarrow \infty }%
\mathbb{E}_{\mathbf{x},a}^{1}\left[ Y_{k}|\mu >\gamma n\right] =\widehat{%
\mathbb{E}}_{\mathbf{x},a}^{1}\left[ Y_{\infty }\right]
\end{equation*}%
for any $\gamma >1$. Hence, writing for $k<n$ and $\gamma >1$%
\begin{equation*}
\mathbb{E}_{\mathbf{x},a}^{1}\left[ Y_{n}|\mu >n\right] =\mathbb{E}_{\mathbf{%
x},a}^{1}\left[ Y_{n}|\mu >\gamma n\right] \frac{\mathbb{P}_{\mathbf{x}%
,a}^{1}\left( \mu >\gamma n\right) }{\mathbb{P}_{\mathbf{x},a}^{1}\left( \mu
>n\right) }+\frac{\mathbb{E}_{\mathbf{x},a}^{1}\left[ Y_{n};\mu >n,\mu \leq
\gamma n\right] }{\mathbb{P}_{\mathbf{x},a}^{1}\left( \mu >n\right) },
\end{equation*}%
observing that
\begin{equation*}
\frac{\mathbb{E}_{\mathbf{x},a}^{1}\left[ Y_{n};\mu >n,\mu \leq \gamma n%
\right] }{\mathbb{P}_{\mathbf{x},a}^{1}\left( \mu >n\right) }=O\left( \frac{%
\mathbb{P}_{\mathbf{x},a}^{1}\left( \mu >n,\mu \leq \gamma n\right) }{%
\mathbb{P}_{\mathbf{x},a}^{1}\left( \mu >n\right) }\right)
\end{equation*}%
and using Lemma~\ref{L_impor} we conclude that, as $n\rightarrow \infty $
\begin{equation*}
\mathbb{E}_{\mathbf{x},a}^{1}\left[ Y_{n}|\mu >n\right] =\left( \widehat{%
\mathbb{E}}_{\mathbf{x},a}^{1}\left[ Y_{\infty }\right] +o(1)\right) \left(
\frac{1}{\sqrt{\gamma }}+o(1)\right) +O\left( 1-\frac{1}{\sqrt{\gamma }}%
\right) .
\end{equation*}%
Letting now sequentially $n$ to infinity and $\gamma $ to $1$ completes the
proof of the lemma.

The next lemma is a generalization of Lemma 3.1 of \cite{LPP2016} to our
setting.

\begin{lemma}
\label{L_Fourheadd1} Under the conditions of Theorem \ref{thm:nonextinct}
for any pair $\left( \mathbf{x},a\right) \in \mathbb{X}\times (-\infty ,0)$
\begin{equation*}
\widehat{\mathbb{E}}_{\mathbf{x},a}^{1}\left[ \sum_{n=1}^{\infty }\mathcal{T}%
_{n}e^{S_{n}}\right] <\infty .
\end{equation*}
\end{lemma}

The proof of this lemma has practically no differences with the proof of
Lemma 4 in \cite{VD2017} and we omit it.

\section{Proof of Theorem \protect\ref{thm:nonextinct}}

For every environmental sequence $\mathbf{F}_{n}$ and $0\leq k<n$ define
\begin{eqnarray*}
\mathbf{F}_{k,n}(\mathbf{s}) &=&\left( F_{k,n}^{(1)}(\mathbf{s}%
),...,F_{k,n}^{(p)}(\mathbf{s})\right) :=\mathbf{F}_{k+1}(\mathbf{F}%
_{k+2}(\ldots (\mathbf{F}_{n}(\mathbf{s}))...)), \\
\mathbf{F}_{n,k}(\mathbf{s}) &=&\left( F_{n,k}^{(1)}(\mathbf{s}%
),...,F_{n,k}^{(p)}(\mathbf{s})\right) :=\mathbf{F}_{n}(\mathbf{F}%
_{n-1}\ldots (\mathbf{F}_{k+1}(\mathbf{s}))...))
\end{eqnarray*}%
and set $\mathbf{F}_{n,n}(\mathbf{s}):=\mathbf{s}$. It is immediate from the
definition of the process $\mathcal{Z}$ that
\begin{equation*}
\mathbb{E}[\mathbf{s}^{\mathbf{Z}_{n}}|\mathbf{Z}_{0}=\mathbf{e}_{i}]=%
\mathbb{E}[F_{0,n}^{(i)}(\mathbf{s})].
\end{equation*}%
Letting $\mathbf{s}=\mathbf{0}$ and using the independency of the
environmental components we get
\begin{equation*}
\mathbb{P}\left( |\mathbf{Z}_{n}|>0\big|\mathbf{Z}_{0}=\mathbf{e}_{i}\right)
=1-\mathbb{E}[F_{0,n}^{(i)}(\mathbf{0})]=\mathbb{E}[1-F_{n,0}^{(i)}(\mathbf{0%
})].
\end{equation*}

Set%
\begin{equation*}
\mathbf{L}_{n,k}:=\mathbf{M}_{n}\mathbf{M}_{n-1}\ldots \mathbf{M}_{k},\
1\leq k\leq n
\end{equation*}%
and let $\mathbf{L}_{n,n+1}$ be the $p\times p$ identity matrix.

We take $\theta =1$ in (\ref{p-def}) and apply the corresponding change of
measure to the representation
\begin{equation*}
1-F_{0,n}^{(i)}(\mathbf{\mathbf{0}})\overset{d}{=}1-F_{n,0}^{(i)}(\mathbf{%
\mathbf{0}})=(\mathbf{e}_{i},\mathbf{1}-\mathbf{F}_{n,0}(\mathbf{\mathbf{0}}%
)).
\end{equation*}%
From now on we agree to consider $\mathbf{e}_{i}$\ as a row vector, and $%
\mathbf{e}_{i}^{T}$ as its transpose. Since $(\mathbf{e}_{i},\mathbf{1}-%
\mathbf{F}_{n,0}(\mathbf{\mathbf{0}}))$ is measurable with respect to $%
\mathcal{F}_{n}$, it follows that
\begin{align*}
\mathbb{E}[1-F_{0,n}^{(i)}(\mathbf{\mathbf{0}})]& =\lambda ^{n}(1)r_{1}(%
\mathbf{e}_{i}^{T})\mathbb{E}\left[ p_{n}^{1}(\mathbf{e}_{i}^{T},\mathbf{L}%
_{n,1})\frac{(\mathbf{e}_{i},\mathbf{1}-\mathbf{F}_{n,0}(\mathbf{\mathbf{0}}%
))}{|\mathbf{L}_{n,1}\mathbf{e}_{i}^{T}|r_{1}(\mathbf{L}_{n,1}\cdot \mathbf{e%
}_{i}^{T})}\right] \\
& =\lambda ^{n}(1)r_{1}(\mathbf{e}_{i}^{T})\mathbb{E}^{1}\left[ \frac{(%
\mathbf{e}_{i},\mathbf{1}-\mathbf{F}_{n,0}(\mathbf{\mathbf{0}}))}{|\mathbf{L}%
_{n,1}\mathbf{e}_{i}^{T}|r_{1}(\mathbf{L}_{n,1}\cdot \mathbf{e}_{i}^{T})}%
\right] .
\end{align*}

To prove \ the theorem we need to show that there exist positive constants $%
C^{-}$ and $C^{+}$ such that%
\begin{equation*}
\frac{C^{-}}{\sqrt{n}}\leq \mathbb{E}^{1}\left[ \frac{(\mathbf{e}_{i},%
\mathbf{1}-\mathbf{F}_{n,0}(\mathbf{\mathbf{0}}))}{|\mathbf{L}_{n,1}\mathbf{e%
}_{i}^{T}|r_{1}(\mathbf{L}_{n,1}\cdot \mathbf{e}_{i}^{T})}\right] \leq \frac{%
C^{+}}{\sqrt{n}}
\end{equation*}%
for all $n\geq 1.$

First observe that $r_{1}(\mathbf{x})$ is a positive function on the compact
$\mathbb{X}$. Hence
\begin{equation*}
0<c_{1}\leq r_{1}(\mathbf{x})\leq c_{2}<\infty
\end{equation*}%
for some constants $c_{1}$ and $c_{2}$. Thus, to complete the proof of
Theorem \ref{thm:nonextinct} it is sufficient to demonstrate that
\begin{equation}
\frac{c_{3}}{\sqrt{n}}\leq \mathbb{E}^{1}\left[ \frac{(\mathbf{e}_{i},%
\mathbf{1}-\mathbf{F}_{n,0}(\mathbf{\mathbf{0}}))}{|\mathbf{L}_{n,1}\mathbf{e%
}_{i}^{T}|}\right] \leq \frac{c_{4}}{\sqrt{n}}  \label{Rough_below_above}
\end{equation}%
for some positive constants $c_{3}$ and $c_{4}$.

\textbf{Estimate in (\ref{Rough_below_above}) from above.} We fix a pair $%
\left( \mathbf{x},a\right) \in \mathbb{X}\times (-\infty ,0)$ and use the
decomposition
\begin{eqnarray}
\mathbb{E}^{1}\left[ \frac{(\mathbf{e}_{i},\mathbf{1}-\mathbf{F}_{n,0}(%
\mathbf{\mathbf{0}}))}{|\mathbf{L}_{n,1}\mathbf{e}_{i}^{T}|}\right]  &=&%
\mathbb{E}_{\mathbf{x},a}^{1}\left[ \frac{(\mathbf{e}_{i},\mathbf{1}-\mathbf{%
F}_{n,0}(\mathbf{\mathbf{0}}))}{|\mathbf{L}_{n,1}\mathbf{e}_{i}^{T}|}\right]
\notag \\
&=&\mathbb{E}_{\mathbf{x},a}^{1}\left[ \frac{(\mathbf{e}_{i},\mathbf{1}-%
\mathbf{F}_{n,0}(\mathbf{\mathbf{0}}))}{|\mathbf{L}_{n,1}\mathbf{e}_{i}^{T}|}%
;\mu \leq n\right]   \notag \\
&&+\mathbb{E}_{\mathbf{x},a}^{1}\left[ \frac{(\mathbf{e}_{i},\mathbf{1}-%
\mathbf{F}_{n,0}(\mathbf{\mathbf{0}}))}{|\mathbf{L}_{n,1}\mathbf{e}_{i}^{T}|}%
;\mu >n\right] .  \label{Decompos}
\end{eqnarray}%
Write $\mathbf{L}_{n,k}=\left( l_{n,k}(q,r)\right) _{q,r=1}^{p}$. Note that
if \ Condition $\mathbf{H3}$ is valid then, according to Lemma 2 in \cite%
{FK60} for any $n,k$ and any tuple $1\leq h,g,q,r\leq p$
\begin{equation}
\Delta ^{-2}\leq \frac{l_{n,k}(h,g)}{l_{n,k}(q,r)}\leq \Delta ^{2}.
\label{Kers}
\end{equation}%
Using for $0\leq k\leq n-1$ the inequality
\begin{eqnarray*}
(\mathbf{e}_{i},\mathbf{1}-\mathbf{F}_{n,0}(\mathbf{0})) &\leq &(\mathbf{e}%
_{i},\mathbf{L}_{n,k+1}\left( \mathbf{1}-\mathbf{F}_{k,0}(\mathbf{0})\right)
) \\
&\leq &(\mathbf{e}_{i},\mathbf{L}_{n,k+1}\mathbf{1})\leq
p^{2}\max_{q,r}l_{n,k+1}(q,r)
\end{eqnarray*}%
and the estimate
\begin{equation*}
|\mathbf{L}_{n,1}\mathbf{e}_{i}^{T}|=|\mathbf{L}_{n,k+1}\mathbf{L}_{k,1}%
\mathbf{e}_{i}^{T}|\geq \min_{q,r}l_{n,k+1}(q,r)|\mathbf{L}_{k,1}\mathbf{e}%
_{i}^{T}|\geq \Delta ^{-2}\max_{q,r}l_{n,k+1}(q,r)|\mathbf{L}_{k,1}\mathbf{e}%
_{i}^{T}|
\end{equation*}%
we see that
\begin{equation*}
\frac{(\mathbf{e}_{i},\mathbf{1}-\mathbf{F}_{n,0}(\mathbf{\mathbf{0}}))}{|%
\mathbf{L}_{n,1}\mathbf{e}_{i}^{T}|}\leq p^{2}\frac{\max_{q,r}l_{n,k+1}(q,r)%
}{\min_{q,r}l_{n,k+1}(q,r)}\frac{1}{|\mathbf{L}_{k,1}\mathbf{e}_{i}^{T}|}%
\leq \frac{\Delta ^{2}p^{2}}{|\mathbf{L}_{k,1}\mathbf{e}_{i}^{T}|}.
\end{equation*}%
Hence we deduce
\begin{equation*}
\mathbb{E}_{\mathbf{e}_{i},a}^{1}\left[ \frac{(\mathbf{e}_{i},\mathbf{1}-%
\mathbf{F}_{n,0}(\mathbf{\mathbf{0}}))}{|\mathbf{L}_{n,1}\mathbf{e}_{i}^{T}|}%
;\mu \leq n\right] \leq \Delta ^{2}p^{2}\mathbb{E}_{\mathbf{e}_{i},a}^{1}%
\left[ \min_{0\leq k\leq n}\frac{1}{|\mathbf{L}_{k,1}\mathbf{e}_{i}^{T}|}%
;\mu \leq n\right]
\end{equation*}%
or, in view of $|\mathbf{L}_{k,1}\mathbf{e}_{i}^{T}|=\left\vert \mathbf{e}%
_{i}\left( \mathbf{L}_{k,1}\right) ^{T}\right\vert $
\begin{eqnarray*}
&&\mathbb{E}_{\mathbf{e}_{i},a}^{1}\left[ \frac{(\mathbf{e}_{i},\mathbf{1}-%
\mathbf{F}_{n,0}(\mathbf{\mathbf{0}}))}{|\mathbf{L}_{n,1}\mathbf{e}_{i}^{T}|}%
;\mu \leq n\right]  \\
&&\qquad\leq\Delta ^{2}p^{2}\mathbb{E}_{\mathbf{e}_{i},0}^{1}\left[ \min_{0\leq
k\leq n}\frac{1}{|\mathbf{e}_{i}\left( \mathbf{L}_{k,1}\right) ^{T}|}%
;\max_{0\leq k\leq n}\log \left\vert \mathbf{e}_{i}\left( \mathbf{L}%
_{k,1}\right) ^{T}\right\vert \geq -a\right] .
\end{eqnarray*}%
To evaluate the right-hand side of this inequality we use the estimates
\begin{multline*}
\mathbb{E}_{\mathbf{e}_{i},0}^{1}\left[ e^{-\max_{0\leq k\leq n}\log |%
\mathbf{e}_{i}\left( \mathbf{L}_{k,1}\right) ^{T}|};\max_{0\leq k\leq n}\log
\left\vert \mathbf{e}_{i}\left( \mathbf{L}_{k,1}\right) ^{T}\right\vert \geq
-a\right]  \\
\leq \sum_{j=-a}^{\infty }e^{-j}\mathbb{P}_{\mathbf{e}_{i},0}^{1}\left(
j<\max_{0\leq k\leq n}\log \left\vert \mathbf{e}_{i}\left( \mathbf{L}%
_{k,1}\right) ^{T}\right\vert \leq j+1\right)  \\
\leq \sum_{j=-a}^{\infty }e^{-j}\mathbb{P}_{\mathbf{e}_{i},j}^{1}\left( \mu
>n\right) \leq \frac{C}{\sqrt{n}}\sum_{j=-a}^{\infty }e^{-j}(j+1),
\end{multline*}%
where the last inequality is justified by (\ref{Mu_above}). Whence, for the
first term at the right-hand side of (\ref{Decompos}) we obtain
\begin{equation*}
\mathbb{E}_{\mathbf{e}_{i},a}^{1}\left[ \frac{(\mathbf{e}_{i},\mathbf{1}-%
\mathbf{F}_{n,0}(\mathbf{\mathbf{0}}))}{|\mathbf{L}_{n,1}\mathbf{e}_{i}^{T}|}%
;\mu \leq n\right] \leq \frac{\Delta ^{2}p^{2}C}{\sqrt{n}}%
\sum_{j=-a}^{\infty }e^{-j}(j+1).
\end{equation*}

For the second term in (\ref{Decompos})\ we apply (\ref{Mu_above}) once
again to conclude that
\begin{eqnarray*}
\mathbb{E}_{\mathbf{e}_{i},a}^{1}\left[ \frac{(\mathbf{e}_{i},\mathbf{1}-%
\mathbf{F}_{n,0}(\mathbf{\mathbf{0}}))}{|\mathbf{L}_{n,1}\mathbf{e}_{i}^{T}|}%
;\mu >n\right]  &\leq &\mathbb{E}_{\mathbf{e}_{i},a}^{1}\left[ \frac{(%
\mathbf{e}_{i},\mathbf{L}_{n,1}\mathbf{1})}{|\mathbf{L}_{n,1}\mathbf{e}%
_{i}^{T}|};\mu >n\right]  \\
&\leq &\Delta ^{2}p^{2}\mathbb{P}_{\mathbf{e}_{i},a}^{1}\left( \mu >n\right)
\leq \frac{\Delta ^{2}p^{2}C\left( 1+\left\vert a\right\vert \right) }{\sqrt{%
n}}.
\end{eqnarray*}%
Thus,
\begin{equation*}
\sqrt{n}\mathbb{E}^{1}\left[ \frac{(\mathbf{e}_{i},\mathbf{1}-\mathbf{F}%
_{n,0}(\mathbf{\mathbf{0}}))}{|\mathbf{L}_{n,1}\mathbf{e}_{i}^{T}|}\right]
\leq C\Delta ^{2}p^{2}\left( 1+\left\vert a\right\vert \right) \left(
1+\sum_{j=-a}^{\infty }e^{-j}(j+1)\right)
\end{equation*}%
which leads to the desired estimate from above in (\ref{Rough_below_above}).

\textbf{Estimate in (\ref{Rough_below_above}) from below}. For a generating
function $\mathbf{f}$, the corresponding mean matrix
\begin{equation*}
\mathbf{m}=\left( \frac{\partial f^{(i)}}{\partial s^{j}}(\mathbf{1})\right)
_{i,j=1}^{p},
\end{equation*}%
and a matrix $\mathbf{a}$ with nonnegative elements define
\begin{equation*}
\psi _{\mathbf{f,a}}(\mathbf{s}):=\frac{\left\vert \mathbf{a}\right\vert }{|%
\mathbf{a}\left( \mathbf{1}-\mathbf{f}\left( \mathbf{s}\right) \right) |}-%
\frac{\left\vert \mathbf{a}\right\vert }{|\mathbf{am}\left( \mathbf{1}-%
\mathbf{s}\right) |},\quad \mathbf{s}\in \left[ 0,1\right] ^{p}\backslash
\left\{ \mathbf{1}\right\} .
\end{equation*}

Let $\mathbf{a}_{i}=\left( a(j,k\right) _{j,k=1}^{p}$ be the matrix with $%
a\left( i,i\right) =1$ and $a\left( k,l\right) =0$ for all $(k,l)\neq (i,i)$%
. Then, clearly,
\begin{equation*}
1-F_{n,0}^{(i)}(\mathbf{s})=|\mathbf{a}_{i}(\mathbf{1}-\mathbf{F}_{n,0}(%
\mathbf{s}))|.
\end{equation*}%
Using the definition of $\psi $, we write
\begin{align*}
\frac{1}{1-F_{n,0}^{(i)}(\mathbf{s})}& =\frac{|\mathbf{a}_{i}|}{|\mathbf{a}%
_{i}(\mathbf{1}-\mathbf{F}_{n,0}(\mathbf{s}))|} \\
& =\frac{1}{|\mathbf{a}_{i}\mathbf{M}_{n}(\mathbf{1}-\mathbf{F}_{n-1,0}(%
\mathbf{s}))|}+\psi _{\mathbf{F}_{n},\mathbf{a}_{i}}\left( \mathbf{F}%
_{n-1,0}(\mathbf{s})\right) \\
& =\frac{1}{\left\vert \mathbf{a}_{i}\mathbf{L}_{n,n}\right\vert }\frac{%
\left\vert \mathbf{a}_{i}\mathbf{L}_{n,n}\right\vert }{|\mathbf{a}_{i}%
\mathbf{L}_{n,n}(\mathbf{1}-\mathbf{F}_{n-1,0}(\mathbf{s}))|}+\frac{1}{%
\left\vert \mathbf{a}_{i}\mathbf{L}_{n,n+1}\right\vert }\psi _{\mathbf{F}%
_{n},\mathbf{a}_{i}\mathbf{L}_{n,n+1}}\left( \mathbf{F}_{n-1,0}(\mathbf{s}%
)\right) \\
& =\frac{1}{|\mathbf{a}_{i}\mathbf{L}_{n,n-1}(\mathbf{1}-\mathbf{F}_{n-2,0}(%
\mathbf{s}))|} \\
& +\frac{1}{\left\vert \mathbf{a}_{i}\mathbf{L}_{n,n}\right\vert }\psi _{%
\mathbf{F}_{n-1},\mathbf{a}_{i}\mathbf{L}_{n,n}}\left( \mathbf{F}_{n-2,0}(%
\mathbf{s})\right) +\frac{1}{\left\vert \mathbf{a}_{i}\mathbf{L}%
_{n,n+1}\right\vert }\psi _{\mathbf{F}_{n},\mathbf{a}_{i}\mathbf{L}%
_{n,n+1}}\left( \mathbf{F}_{n-1,0}(\mathbf{s})\right) .
\end{align*}%
Iterating this procedure, we obtain
\begin{equation}
\frac{1}{1-F_{n,0}^{(i)}(\mathbf{s})}=\frac{1}{|\mathbf{a}_{i}\mathbf{L}%
_{n,1}(\mathbf{1}-\mathbf{s})|}+\sum_{k=1}^{n}\frac{1}{|\mathbf{a}_{i}%
\mathbf{L}_{n,k+1}|}\psi _{\mathbf{F}_{k},\mathbf{a}_{i}\mathbf{L}%
_{n,k+1}}\left( \mathbf{F}_{k-1,0}(\mathbf{s})\right) .  \label{RepresDyakon}
\end{equation}

In view of (\ref{RepresDyakon}) we have%
\begin{equation*}
\mathbb{E}_{\mathbf{e}_{i},a}^{1}\left[ \frac{(\mathbf{e}_{i},\mathbf{1}-%
\mathbf{F}_{n,0}(\mathbf{\mathbf{0}}))}{|\mathbf{L}_{n,1}\mathbf{e}_{i}^{T}|}%
;\mu >n\right] =\mathbb{E}_{\mathbf{e}_{i},a}^{1}\left[ \frac{|\mathbf{e}_{i}%
\mathbf{L}_{n,1}\mathbf{1}|}{|\mathbf{L}_{n,1}\mathbf{e}_{i}^{T}|}\Xi
_{n};\mu >n\right] ,
\end{equation*}%
where
\begin{equation*}
\Xi _{n}:=\left( 1+\sum_{k=1}^{n}\frac{|\mathbf{e}_{i}\mathbf{L}_{n,1}%
\mathbf{1}|}{|\mathbf{a}_{i}\mathbf{L}_{n,k+1}|}\psi _{\mathbf{F}_{k},%
\mathbf{a}_{i}\mathbf{L}_{n,k+1}}\left( \mathbf{F}_{k-1,0}(\mathbf{0}%
)\right) \right) ^{-1}.
\end{equation*}

Using (\ref{Kers}) we conclude that
\begin{equation*}
\frac{|\mathbf{e}_{i}\mathbf{L}_{n,1}\mathbf{1}|}{|\mathbf{L}_{n,1}\mathbf{e}%
_{i}^{T}|}\geq \frac{\min_{q,r}l_{n,1}(q,r)}{p\max_{q,r}l_{n,1}(q,r)}\geq
\frac{1}{p\Delta ^{2}}.
\end{equation*}%
Further, it is known (see Lemma 5 in \cite{VD2017}), that, for all $\mathbf{s%
}\in \left[ 0,1\right] ^{p}\backslash \left\{ \mathbf{1}\right\} $
\begin{equation*}
0\leq \psi _{\mathbf{F}_{k},\mathbf{a}_{i}\mathbf{L}_{n,k+1}}\left( \mathbf{s%
}\right) \leq \Delta p^{2}\mathcal{T}_{k}
\end{equation*}%
and, evidently,
\begin{equation*}
\frac{|\mathbf{e}_{i}\mathbf{L}_{n,1}\mathbf{1}|}{|\mathbf{a}_{i}\mathbf{L}%
_{n,k+1}|}=\frac{|\mathbf{e}_{i}\mathbf{L}_{n,k+1}\mathbf{L}_{k,1}\mathbf{1}|%
}{|\mathbf{a}_{i}\mathbf{L}_{n,k+1}|}\leq \frac{\left\vert \mathbf{e}_{i}%
\mathbf{L}_{n,k+1}\right\vert |\mathbf{L}_{k,1}\mathbf{1}|}{|\mathbf{a}_{i}%
\mathbf{L}_{n,k+1}|}=\left\vert \mathbf{L}_{k,1}\right\vert .
\end{equation*}%
Thus, there exists a positive constant $c_{6}$ such that
\begin{eqnarray*}
\Xi _{n} &\geq &Y_{n}:=c_{6}\left( 1+\sum_{k=1}^{n}\left\vert \mathbf{L}%
_{k,1}\right\vert \mathcal{T}_{k}\right) ^{-1} \\
&\geq &c_{6}\left( 1+\sum_{k=1}^{\infty }\left\vert \mathbf{L}%
_{k,1}\right\vert \mathcal{T}_{k}\right) ^{-1}=:Y_{\infty }>0\text{ \ \ \ }%
\mathbb{\hat{P}}_{x,a}^{1}\text{ a.s}.
\end{eqnarray*}%
(the last in view of Lemma \ref{L_Fourheadd1}). Thus,
\begin{eqnarray*}
\mathbb{E}_{\mathbf{e}_{i},a}^{1}\left[ \frac{|\mathbf{e}_{i}\mathbf{L}_{n,1}%
\mathbf{1}|}{|\mathbf{L}_{n,1}\mathbf{e}_{i}^{T}|}\Xi _{n};\mu >n\right]
&\geq &\frac{1}{p\Delta ^{2}}\mathbb{E}_{\mathbf{e}_{i},a}^{1}\left[
Y_{n};\mu >n\right]  \\
&=&\frac{1}{p\Delta ^{2}}\mathbb{E}_{\mathbf{e}_{i},a}^{1}\left[ Y_{n}|\mu >n%
\right] \mathbb{P}_{\mathbf{x},a}^{1}\left( \mu >n\right) .
\end{eqnarray*}%
Hence, using Lemma \ref{L_impor}, Lemma \ref{L_basic} for the sequence $%
Y_{n}\rightarrow Y_{\infty }$ $\widehat{\mathbb{P}}_{\mathbf{x},a}^{1}$ -
a.s. as $n\rightarrow \infty $ and (\ref{HarmonAbove}) we deduce that
\begin{equation*}
\lim \inf_{n\rightarrow \infty }\sqrt{n}\mathbb{E}_{\mathbf{e}_{i},a}^{1}%
\left[ \frac{|\mathbf{e}_{i}\mathbf{L}_{n,1}\mathbf{1}|}{|\mathbf{L}_{n,1}%
\mathbf{e}_{i}^{T}|}\Xi _{n};\mu >n\right] >0,
\end{equation*}%
proving the estimatefrom below in (\ref{Rough_below_above}). This completes
the proof of Theorem~\ref{thm:nonextinct}.


\begin{thebibliography}{99}
\bibitem{4h} \textit{Afanasyev V.I., Geiger J., Kersting G., Vatutin V.A.}
Criticality for branching processes in random environment. -- Ann. Probab.,
33:2 (2005), 645-673.

\bibitem{ABKV2014} \textit{Afanasyev V.I., Boeinghoff Ch., Kersting G.,
Vatutin V.A.} Conditional limit theorems for intermediately subcritical
branching processes in random environment. -- Ann. Inst. H. Poincar$\acute{e}
$ Probab. Statist., 50:2 (2014), 602-627.

\bibitem{BDGM2014} \textit{Buraczewski D., Damek E., Guivarc'h Y., and
Mentemeier S.} On multidimensional Mandelbrot's cascades. -- J.Difference
Equ. Appl., 20:11 (2014), 1523--1567.

\bibitem{CM2016} \textit{Collamore J.F., Mentemeier S.} Large excursions and
conditioned laws for recursive sequences generated by random matrices. --
arXiv: 1608.05175.

\bibitem{Dyak18} \textit{Dyakonova E.E.} A subcritical decomposable
branching process in a mixed environment. -- Discrete Math. Appl., 28:5
(2018), 275--283.

\bibitem{Dyak15} \textit{Dyakonova E.E.} Limit theorem for multitype
critical branching process evolving in random environment. -- Discrete Math.
Appl., 25:3 (2015), 137-147.

\bibitem{Dy2013} \textit{Dyakonova E.E.} Multitype subcritical branching
processes in a random environment. -- Proc. Steklov Inst. Math., 282 (2013),
80--89.

\bibitem{Dyak12} \textit{Dyakonova E. E.} Multitype branching processes
evolving in a Markovian environment. -- Discrete Math. Appl., 22:5-6 (2012),
639--664.

\bibitem{Dyak11} \textit{Dyakonova E.E.} Multitype Galton-Watson branching
processes in Markovian random environment. -- Theory Probab. Appl., 56:3
(2011), 508--517.

\bibitem{VatDyak10} \textit{Vatutin V.A., Dyakonova E.E.} Asymptotic
properties of multitype critical branching processes evolving in a random
environment. -- Discrete Math. Appl., 20:2 (2010), 157--177.

\bibitem{Dyak08} \textit{Dyakonova E.E.} On subcritical multi-type branching
process in random environment. -- Fifth Colloquium on Mathematics and
Computer Science, Discrete Math. Theor. Comput. Sci. Proc., AG, Assoc.
Discrete Math. Theor. Comput. Sci., Nancy, 2008, 397--404.

\bibitem{Dyak07} \textit{Dyakonova E.E.} Critical multitype branching
processes in a random environment. -- Discrete Math. Appl., 17:6 (2007),
587--606.

\bibitem{FK60} \textit{Furstenberg H., Kesten H.} Products of random
matrices. -- Ann. Math. Statist., 31 (1960), 457--469.

\bibitem{KV2017} \textit{Kersting G., Vatutin V.} Discrete Time Branching
Processes in Random Environment, Wiley, John Wiley \& Sons Inc., New Jersey,
USA; ISTE, London, UK, 2017, 306 pp.

\bibitem{LPP2016} \textit{Le Page E., Peigne M., Pham C.} The survival
probability of a critical multi-type branching process in i.i.d. random
environment. -- Ann. Probab., 46:5 (2018), 2946--2972.

\bibitem{Vat2012c} \textit{Vatutin V., Dyakonova E., Jagers P., Sagitov S.}
A decomposable branching process in a Markovian environment. --
International Journal of Stochastic Analysis, 2012, v. 2012, Article ID
694285, 24 p.

\bibitem{VD2017} \textit{Vatutin V., Dyakonova E.\ }Multitype branching
processes in random environment: survival probability for the critical case.
-- Theory Probab. Appl., 62:4 (2018), 506--521.

\bibitem{Vat2013b} \textit{Vatutin V., Liu Q.} Branching processes
evolving in asynchronous environments. -- in: Proceedings 59-th ISI World
Statistics Congress, August 2013, Hong Kong (Hong Kong, 25-30 August 2013),
International Statistical Institute, The Hague, The Netherlands, 2013, 1744--1749.

\bibitem{VZ} \textit{Vatutin V.A., Zubkov A.M.} Branching processes II. --
J. Sov. Math., 67:6 (1993), 3407--3485.

\bibitem{VatWacht2018} \textit{Vatutin V., Wachtel V.} Multi-type
subcritical branching processes in a random environment. -- Adv. in Appl.
Probab., 50:A (2018), 281--289.

\bibitem{PH17} \textit{Pham C.} Conditional limit theorems for products of
positive random matrices. -- arXiv: 1703.04949.
\end{thebibliography}
\end{document}